\documentclass[11pt]{article} 
\usepackage{amsfonts,euscript,theorem}
\newcommand{\no}{\noindent}
\newcommand{\bc}{\begin{center}}
\newcommand{\ec}{\end{center}}
\newcommand{\be}{\begin{equation}}
\newcommand{\ee}{\end{equation}}
\newcommand{\bea}{\begin{eqnarray*}} 
\newcommand{\eea}{\end{eqnarray*}}
\newcommand{\ran}{\ensuremath{\rangle}}
\newcommand{\lan}{\ensuremath{\langle}}
\newcommand{\rar}{\ensuremath{\rightarrow}}
\newcommand{\ii}{\ensuremath{\infty}}

\newcommand{\pdx}{\ensuremath{\partial}_x}
\newcommand{\pdp}{\ensuremath{\partial}_x^{\,p}}
\newcommand{\pdr}{\ensuremath{\partial}_x^{\,}} 
 
\newcommand{\dx}{\ensuremath{\frac{d}{dx}\ }}

\newcommand{\h}[1]{\ensuremath{\hbox{ #1 }}}
 
\newcommand{\w}[1]{\ensuremath{\widetilde{#1}}}

\newcommand{\Om}{\ensuremath{\Omega}}
\newcommand{\si}{\ensuremath{\sigma}}

\newcommand{\ps}{\ensuremath{\psi}}
\newcommand{\ve}{\ensuremath{\varepsilon}}
\newcommand{\vp}{\ensuremath{\varphi}}
\newcommand{\pe}{\ensuremath{\vp _{\ve}}}

 \newcommand{\de}{\ensuremath{\delta}}
\newcommand{\N}{\ensuremath{\mathbb{N}}}
\newcommand{\NN}{\ensuremath{\mathbb{N}_0}}
\newcommand{\R}{\ensuremath{\mathbb{R}}}
\newcommand{\Rm}{\ensuremath{{\mathbb{R}}^m}}
\newcommand{\C}{\ensuremath{\mathbb{C}}}
\newcommand{\A}{\ensuremath{A_0(\R)}}
\newcommand{\Aq}{\ensuremath{A_q(\R)}}
\newcommand{\G}{\ensuremath{{\EuScript{G}}}} 
\newcommand{\GR}{\ensuremath{{\EuScript{G}}(\R)}} 
\newcommand{\Gm}{\ensuremath{{\EuScript{G}}(\Rm)}} 
\newcommand{\D}{\ensuremath{{\EuScript{D}}(\R)}} 
\newcommand{\DD}{\ensuremath{{\EuScript{D}}'(\R)}} 

\newtheorem{L1}{Lemma}
\newtheorem{Th1}{Theorem}
\newtheorem{Th2}[Th1]{Theorem}
\newtheorem{Cor1}{Corollary}
\newtheorem{Cor2}[Cor1]{Corollary}
\newtheorem{Th3}[Th1]{Theorem}

\begin{document}

\sloppy 
\setlength{\baselineskip}{18pt}
\bc 
\boldmath
{\Large \textbf {Balanced Colombeau products \\ of the distributions $x_{\pm}^{-p}$ and $x^{-p}$}}
\unboldmath

\vspace*{3mm}
{\textbf {B. P. Damyanov}}
\ec

\vspace*{8mm}
\setlength{\baselineskip}{15pt}
\no {\small \textbf{Abstract.} \ {Results on singular products of the distributions $x_{\pm}^{-p}$ and $x^{-p}$ \,for natural $p$ \,are derived \,when the products are balanced so that their sum exists in the distribution space. These results follow the pattern of a known distributional product published by Jan Mikusi\'nski in 1966.  The results are obtained in Colombeau algebra of generalized functions, which is most relevant algebraic construction for tackling nonlinear problems of Schwartz distributions.

\vspace*{2mm}
\no \textbf{Keywords:} \ \textit{Schwartz distributions, multiplication, Colombeau generalized functions}

\vspace*{1mm}
\no \textbf{AMS subject classification:} \ 46F10, 46F30}}

\vspace*{5mm}
\setlength{\baselineskip}{17pt}
\textsc{0. INTRODUCTION}

The Colombeau algebra of generalized functions \G \,, introduced first in \cite{colm}, has followed various constructions of differential algebras aimed at solving certain nonlinear problems of Schwartz distributions. An important reason for the growing popularity of the algebra \G \ is its almost optimal properties regarding the long-standing problem of multiplication of distributions. Indeed, \G \ is an associative differential algebra, the multiplication is compatible with products of  $C^\ii$-differentiable functions, and the linear embedding of the distribution space commutes with partial differentiation. Moreover, the so-called 'association' in \G \,, which is a faithful generalization of the equality of distributions, enables obtaining results `on distributional level'.  

In 1966, Jan Mikusi\'nski published in \cite{mik} his well-known result\,: 
\be x^{-1}\,.\,x^{-1} \,-\, \pi^{{}2}\,\de(x)\,.\, \de(x) \ = \ x^{-2}, \ \ x \in \R.  \label{mik} \ee
\no Although neither of the products on the left-hand side here exists, their difference still has a correct meaning in the distribution space \DD. \ Formulas including balanced products of distributions with coinciding singularities can be found in the mathematical and physical literature. We termed such equations `products of Mikusi\'nski type' in a previous paper  \cite{ijpm}, where we derived a generalization of (\ref{mik}) in Colombeau algebra of tempered generalized functions \,(see equation (\ref{fstep}) below). 

Following the pattern of the basic Mikusi\'nski product (\ref{mik}), we prove in this paper further results on balanced  products of the distributions with singular point support $x_{\pm}^{-p}, x^{-p}$, and $\de^{(p)}(x), \,p\in\NN$ \,and $x\in\R$. We evaluate the products as the distributions are embedded in Colombeau algebra and prove that each of the products admits an associated distribution. 

\vspace*{4mm}
\bc \textsc{1. NOTATION AND DEFINITIONS} \ec 

\vspace*{2mm}
\textbf{1.1.} We start with recalling the fundamentals of Colombeau algebra on the real line.

\textit{Notation 1.} If  \N \ stands for the natural numbers, denote $\NN = \N \cup \{0\}$ \,and $\de _{ij} = \{\,1$ \,if \,$i=j, \, = 0$ otherwise$\,\}, \,i, j\in\NN$. Then we put for arbitrary $q \in \NN$\,: 
\[\Aq = \{ \vp(x) \in \D: \int_{\R} x^{j}\,\vp (x)\,dx  = \de _{0j}, \ j = 0, 1,...,q \}.\]  
Set also \,$\pe = \ve ^{-1}\vp (\ve ^{-1}x)$ \,for  $\vp \in \Aq, \ve > 0$ \,and $\check{g}(x) = g(-x)$. Note finally that the shorthand notation $\pdx = \dx$ \,will be used in the one-dimensional case too. 

\vspace*{1mm}
\textit{Definition 1.} Let $\EuScript{E}\,[\R]$ be the algebra of functions $f(\vp , x): \A \times \R  \rar \C$ that are infinitely differentiable for  fixed 'parameter' \vp. Then, the generalized functions of Colombeau are elements of the quotient algebra 
\[ \G \equiv\GR = \EuScript{E}_{\mathrm{M}} [\R]\,/ \ \EuScript{I}\,[\R].\] 
Here $\EuScript{E}_{\mathrm{M}}[\R]$ \,is the subalgebra of `moderate' functions such that for each compact subset $K$ of \R \,and $p \in\N$ \,there is a $q\in \N$ \,such that, for each $\vp \in \Aq$,  
\[ \sup_{x \in K}\,|\pdp \,f(\pe, x)\,| = O(\ve^{-q}), \hbox{ as}  \ \ve \rar 0_+.\] 
The ideal $\EuScript{I}\,[\R]$  of $\EuScript{E}_{\mathrm{M}}[\R]$ consists of all functions such that for each compact $K \subset\R $ \ and any $p\in \N$ \,there is a $q\in\N$ \ such that, for every $r \geq q$ and $\vp \in A_r(\R)$, 
\[ \sup_{x \in K}\,| \pdp \,f(\pe, x)\,| = O(\ve^{r-q}), \hbox{ as} \ \ve \rar 0_+.\] 

The algebra \GR \ contains the distributions on \R, canonically embedded as a \C-vector subspace \,by the map  \vspace*{-1mm}
\be i : \DD \rar \,\G : u \mapsto \w{u} = \{\,\w{u}(\vp, x) := (u*{\check{\vp}})(x): \,\vp\in\Aq \,\}. \label{emb}\ee
The derivative in \GR \ is in consistency with this embedding of distributions\,:
\be \pdx \,\w{u} \ = \ \w{\pdx u},  \quad u\in \DD. \label{cons} \ee 

The equality of generalized functions in  \G \ is very strict and a weaker form of equality {\em{in the sense of association}} is introduced, which plays a fundamental role in Colombeau theory.

\vspace*{1mm} 
\no \textbf{Definition 2.} A generalized function $f\in\GR$ \ is said to be `associated' with (a) another function $g\in\G$, denoted $f \approx g$, \,or (b) a distribution $ u\in \DD$  ($f \approx u$) \ if for some representatives $f(\vp_\ve, x), g(\vp_\ve, x)$ \,and arbitrary $\ps (x)\in\D$ \,there is a $q\in\NN$ \,such that, for any $\vp (x)\in\Aq$, \ it holds $\lim_{{}\ve \rar 0_+} \int_{\R} [f(\vp_\ve , x)- g(\vp_\ve, x)] \ps (x)\,dx = 0$, or respectively $\lim_{{}\ve \rar 0_+} \int_{\R} f(\vp_\ve , x) \ps (x)\,dx = \lan u, \ps\ran.$ 

\vspace*{1mm}
These definitions are independent of the representatives chosen, and the association is a faithful generalization of the equality of distributions \cite{col84}; which implies the following equivalence relation for the embedding of distributions\,: 
\be f \approx \w{u} \Longleftrightarrow f \approx u \quad \quad \h{for each} f \in \GR, \ u \in\DD. \label{equiv}\ee

Now, by Colombeau product of two distributions is meant the product of their embeddings in \GR \ whenever the result admits an associated distribution. 

The following coherence result holds \cite[Proposition 10.3]{mob} : If the regularized model product (in the terminology of Kami\'nski) of two distributions exists, then their Colombeau product also exists and coincides with the former. Moreover, in the general setting of Colombeau algebra \Gm \ \cite{col84} \ (when the parameter functions \vp \ are not defined as tensor products), as well as in the algebra \GR \ on the real line,  this assertion turns into an equivalence, according to a result by Jel\'{\i}nek \cite{jel}; cf. also the recent study on Colombeau algebra in \cite{chi}.

Denote now by $\w{x^{\,-p}}$ and $\w{\de}^{(p)} (x) $ the embeddings (\ref{emb}) in \G \ of the distributions $x^{\,-p}$ and $\de^{(p)}(x), \,p \in \N $. Then the following balanced distributional product in Colombeau algebra was proved in  \cite{ijpm}, which generalizes the basic Mikusi\'{n}ski formula (\ref{mik}) for arbitrary $p,q \in \N$:
\begin{equation} \quad \w{ x^{\,-p}}\,.\,\w{ x^{\,-q}} - \pi^2\frac{(-1)^{p+q}}{(p-1)!\,(q-1)!}\,\w{\de}^{\,(p-1)}(x)\,.\,\w{\de}^{\,(q-1)}(x) \ \approx \ x^{\,-p-q} \ (x\in \R). \label{fstep}\end{equation} 

\vspace*{2mm}
\textbf{1.2.} Consider next the definition of the distributions in study.

\vspace*{1mm}
\textit{Notation 2.} If $a\in\C$ and Re~$a > -1$, denote as usual the locally-integrable functions\,:
\[ x_+^{\,a} =  \{ x^{\,a} \ \h{if} x > 0 , \quad  = 0 \ \h{if} x<0\}, \qquad x_-^{\,a} =  \{  (-x)^{\,a} \ \h{if} x < 0,  \quad = 0 \ \h{if} x >0 \}.  \]
\[ \ln x_+ =  \{ \ln x \ \h{if} x > 0 , \quad  = 0 \ \h{if} x<0\}, \quad \ln x_- =  \{ \ln (-x) \ \h{if} x < 0,  \quad = 0 \ \h{if} x >0 \}.  \]
\[  \ln |x| = \ln x_+ + \ln x_-, \qquad \ln |x|\,\mathrm{sgn}\,x = \ln x_+ - \ln x_-.\]

\vspace*{1mm}
The distributions $x_{\pm}^{\,a}$ \,are defined for any $a\in\Om := \{a \in \R: \, a \ne -1, -2,\ldots\}$, by setting 
\[ x_{+}^{\,a} = \pdr \,x_{+}^{\,a +r }(x), \qquad  x_{-}^{\,a} = (-1)^r \,\pdr \,x_{-}^{\,a +r }(x), \] 
where $r\in \NN$  \,is such that $a+r > -1$ \,and the derivatives are in distributional sense.

This definition can be extended also for negative integer values of $a$ by a procedure essentially due to M. Riesz (see \cite[\S \,3.2]{hor}). For each $\ps(x) \in \D, \ a\mapsto\lan x_{+}^a, \ps\ran$ is an analytic function of $a$ on the set $\Om$. The excluded points are simple poles of this function. For any $p\in\NN$, the residue at $a = -p-1$ \,is \ $\lim_{a\rightarrow -p-1}(a+p+1) \ \lan x_{+}^a, \ps\ran = \psi^{(p)}(0)/p!$. Subtracting the singular part, one gets for any $p\in\NN$\,:
\[\lim_{a\rightarrow -p-1} \ \lan x_{+}^a, \ \ps\ran - \frac{1}{p!}\ \psi^{(p)}(0) = - \frac{1}{p!}\int_0^\infty \ln x \,\ps^{(p)}\,dx + \frac{\psi^{(p)}(0)}{p!} \sum_{k=1}^{p} \frac{1}{\,k}.\]
The right-hand side of this equation, which is the principal part of the Laurent expansion, was proposed by H\"{o}rmander in \cite{hor} to define the distribution $x_+^{-p-1}$, acting here on the test-function $\ps(x)$. In view of Notation~2, this is equivalent to the following definition of $x_{+}^{-p-1}$ \,for arbitrary $p\in\NN \ (x\in\R)$\,:
\[ x_{+}^{-p-1} =  \frac{(-1)^p}{p!} \ \pdx^{p+1} \ln x_+ \ + \ \frac{(-1)^p\,\si_p}{p!} \ \de^{(p)}(x).\]
We have introduced here the shorthand notation
\be \si_p \ := \ \sum_{k=1}^p 1/k  \quad ( p\in\NN), \qquad \hbox{noting that} \ \si_0 = 0.\label{sig}\ee 
Similar arguments lead to the defining equation 
\be x_{-}^{-p-1} =  \frac{-1}{p!} \ \pdx^{p+1} \ln x_- \ + \ \frac{\si_p}{p!} \ \de^{(p)}(x). \label{x-p}\ee
Note that this definition exactly coincides with that of the distributions $x_{\pm}^{\,-p-1} \equiv F_{-p-1}(x_{\pm}, \lambda)$ \,introduced by Gelfand and Shilov, as regularization of the integrals $\int_{\R_{\pm}} x^{\lambda}\,\ps(x)\,dx$ \,taken at the points $\lambda = - p-1$ \ \cite[\S \,1.4]{gsh}. 

\vspace*{1mm}
\no One checks that the distributions $x_{\pm}^{-p}$ satisfy 
\be \pdx \,x_+^{-p} \ = \ -p \ x_+^{-p-1} \ + \ \frac{(-1)^{p}}{p!} \ \de^{(p)}(x), \quad \pdx \,x_-^{-p} \ = \ p \ x_-^{-p-1} \ - \ \frac{1}{p!} \ \de^{(p)}(x). \label{dx-p}\ee
Now, it follows immediately that
\be   x_{+}^{-p}\,|_{x \mapsto -x} = x_{-}^{-p} \qquad \hbox{and} \qquad x_{+}^{-p} \ + \ (-1)^{p} \,x_{-}^{-p} \ = \  x^{-p}. \label{x+-p}\ee
Here the distribution $x^{-p}$ \,is defined, as usual, as distributional derivative of order $p$\,:
\be x^{-p} = \frac{(-1)^{p-1}}{(p-1)!} \ \pdp \ln |x|, \qquad \hbox{and it holds}\quad \pdx x^{-p} \ = \ -p \ x^{-p-1}. \label{defx-p} \ee
For later use we note that the following basic property of the ditrsibutions $x_{\pm}^a$ is preserved for arbitrary $a\in \C$ \,\cite[\S \,3.2]{hor} :
\be x\,.\, x_{\pm}^a \ = \ x_{\pm}^{a+1}. \label{xx} \ee

Recall finally the definition of the distributions $(x {\pm}i 0)^{-p-1}$ for $p\in \NN$ \,and $\,\ x\in\R$\,:
\be (x {\pm}i 0)^{-p-1} := \lim_{y\rar 0_+}  (x {\pm}i y)^{-p-1} = x^{\,-p-1} \ \mp \,\frac{( - 1)^p \ i \,\pi}{p!} \ \de^{(p)}(x).   \label{xip}\ee

\vspace*{3mm}
\bc  \textsc{2. PRELIMINARY RESULTS} \ec 

\vspace*{2mm}
We first recall two results on distributional products in Colombeau algebra that will be used later in the work. Let $\w{x^{\,-p}}, \w{x_{+}^p}$, and $\w{\de}^{(p-1)} (x) $ denote the embeddings in \GR \ of the distributions $x^{\,-p}, x_{+}^p$, and $\de^{(p-1)} (x), \,p \in \N$. Then the following 'ordinary' Colombeau product, given here in dimension one, was obtained in \cite{cmuc}\,:
\be \w{x_{+}^p}\,.\,\w{\de}^{(p)}(x) \ \approx \ ( - 1 )^{\,p} \ \frac{p\,!}{2} \ \de(x) \qquad (x\in\R). \label{cmuc}\ee

\no Further, this balanced Colombeau product was proved in  \cite{ijpm}\,: \,For arbitrary $p,q \in \N$, 
\be \frac{(-1)^{q-1}}{(q-1)!} \ \w{x^{-p}}\,.\,\w{\de}^{(q-1)}(x) +  \frac{(-1)^{p-1}}{(p-1)!} \ \w{x^{- q}}\,.\,\w{\de}^{(p-1)}(x) \approx \frac{(-1)^{p+q-1}}{(p+q-1)!} \ \de^{(p+q-1)}(x). \label{xpdq}\ee

\vspace*{1mm}
Next, we prove a general property of balanced products in Colombeau algebra that will be needed in the sequel too.

\begin{L1}$\!\!.$ The derivative of any balanced Colombeau product of distributions $\sum_{k=1}^{2} \left(\w{u_k}\,.\,\w{v_k}\right) \,\approx \,w \ (u_k, v_k, w \in \D')$ \ admits also an associated distribution and it holds
\be \sum_{k=1}^{2} \left(\w{\pdx u_k}\,.\,\w{v_k}\, + \,\w{u_k}\,.\,\w{\pdx v_k}\right) \,\approx \,\pdx w. \label{l1}\ee
\end{L1}

\textsc{Proof}\,: For a given $\vp \in \A$ \,the representatives $\w{u_k}(\pe, x), \w{v_k}(\pe, x)$ are smooth functions of $x$, by a fixed $\ve$. Therefore, choosing an arbitrary $\ps(x) \in \D$ \,and taking into account equation (\ref{cons}), applied to the representatives of the embeddings, we obtain
\begin{eqnarray} I & := & \int_{- \infty}^{\infty} \ps(x)\,\pdx \!\left( \sum_{k=1}^{2} \left[\,\w{u_k}(\pe, x) \w{v_k}(\pe, x)\right]\right)\,dx \nonumber \\ & = &  \int_{- \infty}^{\infty} \ps(x)\,\sum_{k=1}^{2}\left[\,\w{\,\pdx u_k}(\pe, x)\,\w{v_k}(\pe, x) + \w{u_k}(\pe, x)\,\w{\,\pdx v_k}(\pe, x)\right]\,dx.\label{I}\end{eqnarray}

\no On the other hand, we get on integration by parts\,:
\begin{eqnarray} I & = & \sum_{k=1}^{2} \ \int_{- \infty}^{\infty} \ps(x)\,\pdx \left[\,\w{u_k}(\pe, x)\,\w{v_k}(\pe, x)\right]\,dx \nonumber \\ & = & \sum_{k=1}^{2} \left(\w{u_k}(\pe, x)\,\w{v_k}(\pe, x)\,\ps(x)\,|_{-\infty}^{\ \,\infty} \ - \  \int_{- \infty}^{\infty} \w{u_k}(\pe, x)\,\w{v_k}(\pe, x)\,\ps'(x)\,dx \right) \nonumber \\ & = & \ - \  \int_{- \infty}^{\infty} \sum_{k=1}^{2} \left[\w{u_k}(\pe, x)\,\w{v_k}(\pe, x)\right]\,\pdx \ps(x)\,dx. \label{I'}\end{eqnarray}

\no From equations (\ref{I}), (\ref{I'}) and by the assumption of the theorem, it now follows 
\bea \lim_{{}\ve \rar 0_+} I &=& \lim_{{}\ve \rar 0_+}  \int_{- \infty}^{\infty} \ps(x)\,\sum_{k=1}^{2}\left[\,\w{\,\pdx u_k}(\pe, x)\,\w{v_k}(\pe, x) + \w{u_k}(\pe, x)\,\w{\,\pdx v_k}(\pe, x)\right]\,dx \\ & = & - \lim_{{}\ve \rar 0_+}  \int_{- \infty}^{\infty} \sum_{k=1}^{2} \left[\w{u_k}(\pe, x)\,\w{v_k}(\pe, x)\right]\,\pdx \ps(x)\,dx  =  \ - \ \lan w, \,\pdx\ps\ran \ = \ \lan \pdx w, \,\ps\ran. \eea
According to Definition~2, this proves the existence of an associated distribution for the derivative, as well as equation (\ref{l1}).

\vspace*{3mm}
\bc \textsc{2. MAIN RESULTS} \ec 

\vspace*{2mm}
We now proceed to particular balanced products of distributions obtained in Colombeau algebra \GR. With the notation (\ref{sig}), the following assertion holds.
\begin{Th1}$\!\!.$ For each $p\in\NN$, the embeddings in \GR \,of the distributions $x_{\pm}^{-p}, x_{\pm}^{\,p}$, and $\de(x)$ satisfy\,:
\be (-1)^p \ \w{x_{-}^{\,-p-1}}\,.\,\w{x_{+}^p} \ - \ \w{\ln x_{-}}\,.\,\w{\de}(x) \  \approx \ \frac{\si_p}{2} \ \de(x).   \label{th1-} \ee
\be (-1)^p \ \w{x_{+}^{\,-p-1}}\,.\,\w{x_{-}^p} \ - \ \w{\ln x_{+}}\,.\,\w{\de}(x) \  \approx \ \frac{\si_p}{2} \ \de(x).   \label{th1+} \ee
\end{Th1}

\no \textsc{Proof}\,:   (i) For given $\vp \in \A$, suppose without lost of generality that \,supp $\!\vp(x) \subseteq [-l, l]$ \,for some $l\in \R_+$. Then, the embedding rule (\ref{emb}) and the substitution $u = (y-x)/{\ve}$ give for the representatives
\be  \w{x_+^p}(\pe, x) = \ve ^{-1} \int_0^{-\ve l +x} y^p\,\vp ((y-x)/\ve)\,dy = \int_{- x/ \ve}^{\,l} (\ve u + x)^p\ \vp(u)\,du, \label{wx+}\ee
and 
\be \w{\de}^{(p)}(\pe, x) = \frac{(-1)^p}{\ve^{p+1}}\ \vp^{(p)}\left( - \frac{x}{\ve}\right).\label{hde}\ee
Similarly, equation (\ref{x-p}), the rules (\ref{emb}) for embedding and (\ref{cons}) for Colombeau derivative, as well as the substitution $v = (y-x)/{\ve}$ \,yield
\be \w{x_{-}^{-p-1}}(\pe, x)  = \frac{(-1)^p}{p!\,\ve^{p+1}} \int_{- l}^{- x/\ve} \ln (- \ve v - x)\,\vp^{(p+1)}(v)\,dv  + \frac{(-1)^p\ \si_p}{p!\,\ve^{p+1}} \,\vp^{(p)}\left( - \frac{x}{\ve}\right).\label{wx-p} \ee

For a given $\ps(x) \in \D$, \,evaluate now
\[ F_p  :=  (-1)^{ p} \ \int_{- \infty}^{\infty} \ps(x)\,\w{x_{-}^{-p-1}}(\pe, x)\,\w{x_+^p}(\pe, x)\,dx. \]

\no Inserting equations (\ref{wx+}) and (\ref{wx-p}), we obtain
\bea  F_p & =&  \frac{1}{p!\,\ve^{p+1}} \int_{- \ve l }^{\ve l } dx\,\ps(x) \int_{- x/ \ve}^{\,l} du\ (\ve u + x)^p\ \vp(u) \int_{- l}^{- x/ \ve} \ln (-\ve v - x)\ \vp^{(p+1)}(v)\,dv \\ & & +  \frac{\si_p}{p!\,\ve^{p+1}} \int_{- \ve l }^{\ve l } dx\,\ps(x) \vp^{(p)}\left( - \frac{x}{\ve}\right) \int_{- x/ \ve}^{\,l} (\ve u + x)^p\ \vp(u)\,du =:  \ F_p{\,'} \ + \ F_p{\,''}. \eea
Applying the substitution $w=-\,x/\ve$, Taylor theorem, change of the order of integration, and finally the  substitution $w \,\rar \,t = (w-v)/(u-v)$ \,we get for the first term
\bea F_p{\,'} & = & \frac{1}{p!} \int_{- l}^{\,l} dw\,\ps(-\,\ve w) \int_{\,w}^{\,l} du\ (u-w)^p \vp(u) \int_{- l}^{\,w} \ln (\ve w - \ve v)\,\vp^{(p+1)}(v)\,dv\\ & = & \frac{\ps(0)}{p!} \int_{- l}^{\,l} du\,\vp(u) \int_{-l}^{\,u}dv\,\vp^{(p+1)}(v) \int_{\,v}^{\,u} \ (u-w)^p \ln (\ve w - \ve v)\,dw  + O\,(\ve) \\ & = &\frac{\ps(0)}{p!} \int_{- l}^{\,l} du\,\vp(u) \int_{-l}^{\,u}dv\,\vp^{(p+1)}(v) (u-v)^{p+1} \\ & & \times\left[ \ln (\ve u - \ve v) \int_0^1(1-t)^p\,dt + \int_0^1(1-t)^p\,\ln t \ dt \right] +  O\,(\ve). \eea
Now, we have 
\[ \int_0^1(1-t)^p\,dt = \frac{1}{p+1}, \ \ \hbox{and} \ \int_0^1(1-t)^p\,\ln t \ dt = - \,\frac{\sigma_{p+1}}{p+1}\]
(cf. \cite[\S\,4.24.3]{gr} for the calculation of the second integral). Therefore,
\[ F_p{\,'} = \frac{\ps(0)}{(p+1)!} \int_{- l}^{\,l} du\,\vp(u) \int_{-l}^{\,u} \vp^{(p+1)}(v) (u-v)^{p+1}\left[ \ln (\ve u - \ve v) - \sigma_{p+1}\right] \,dv  +  O\,(\ve).\]
\no Further, the substitution $v=-\,x/\ve$, Taylor theorem, change of the order of integration, and integration by parts in the variable $v$ (the integrated part being 0) give
\bea  F_p{\,''} & = & \frac{\sigma_p}{p!} \int_{- l}^{\,l} dv\,\ps(-\,\ve v) \vp^{(p)}(v) \int_{\,v}^{\,l} (u-v)^p \vp(u)\,du \\ & = & \frac{\sigma_p \ \ps(0)}{p!} \int_{- l}^{\,l} du\,\vp(u) \int_{-l}^{\,u} (u-v)^p \,\vp^{(p)}(v)\,dv + O\,(\ve)  \\ & = & \frac{\sigma_p \ \ps(0)}{(p+1)!} \int_{- l}^{\,l} du\,\vp(u) \int_{-l}^{\,u} (u-v)^{p+1} \,\vp^{(p+1)}(v)\,dv + O\,(\ve). \eea
Note that to obtain the asymptotic evaluations of $F_p{\,'}$ and $F_p{\,''}$, we have taken into account that the second  term in each Taylor expansion is multiplied by definite integrals majorizable by constants. Replacing the obtained expressions into $F_p$ and integrating by parts, we get
\begin{eqnarray}  F_p & = & \frac{\ps(0)}{(p+1)!} \int_{- l}^{\,l} du\,\vp(u) \int_{-l}^{\,u} \vp^{(p+1)}(v) (u-v)^{p+1}\left[ \ln (\ve u - \ve v) - \frac{1}{p+1}\right] dv  +  O\,(\ve)\nonumber \\ & = &  \frac{\ps(0)}{p!} \int_{- l}^{\,l} du\,\vp(u) \int_{-l}^{\,u}\vp^{(p)}(v) (u-v)^{p} \,\ln (\ve u - \ve v)\ dv  +  O\,(\ve). \label{F_p}\end{eqnarray}
Denoting further by $I_{p}$ \,the second integral in this equation divided by $p!$, we get on integration by parts 
\[ I_p = \frac{1}{(p-1)!} \int_{-l}^{\,u}\vp^{(p-1)}(v) (u-v)^{p-1} \,\ln (\ve u - \ve v)\ dv + \frac{1}{p!} \int_{-l}^{\,u}dv\,\vp^{(p-1)}(v) (u-v)^{p-1}.\]
Iterating this procedure $p$ times, taking into account that for each $p\in\N$ \ it holds
\[ \frac{1}{p!} \int_{-l}^{\,u} \vp^{(p)}(v) (u-v)^{p}\,dv = \int_{-l}^{\,u} \vp(v)\,dv, \]
we obtain
\[ I_p = \int_{-l}^{\,u} \,\ln (\ve u - \ve v)\ dv + \sigma_p \int_{-l}^{\,u} \vp(v)\,dv.\]
The replacement of $I_p$ in equation (\ref{F_p}) now gives
\[F_p =  \ps(0) \int_{- l}^{\,l} du\,\vp(u) \int_{-l}^{\,u} \ln (\ve u - \ve v)\,\vp(v)\,dv \ + \ \sigma_p \,\ps(0) \int_{- l}^{\,l} du\,\vp(u) \int_{-l}^{\,u} \vp(v)\,dv \ + \ O\,(\ve). \]  

(ii) On the other hand, equation (\ref{hde}) in the case $p=0$, the substitution $u= - x/\ve$, and Taylor theorem yield
\bea G & := & \int_{- \infty}^{\infty} \ps(x)\,\w{\ln x_-}(\pe, x)\,\w{\de}(\pe, x)\,dx \nonumber \\ &=&  \int_{- l}^{\,l} du\,\ps(-\,\ve u) \,\vp(u) \int_{- l}^{\,u} \ln (\ve u - \ve v)\,\vp(v)\,dv  \\ & = & \ps(0) \int_{- l}^{\,l} du\,\vp(u) \int_{-l}^{\,u} \ln (\ve u - \ve v)\,\vp(v)\,dv \ + \ O\,(\ve).  \eea
Therefore, 
\[F_p \ - \ G = \sigma_p\,\ps(0) \int_{- l}^{\,l} du\,\vp(u) \int_{-l}^{\,u} \vp(v)\,dv = \frac{\sigma_p\ \ps(0)}{2}.\]
By linearity, this implies
\[\lim_{{}\ve \rar 0_+} \int_{- \infty}^{\infty} \ps(x)\,\left[(-1)^p \,\w{x_{-}^{-p-1}}(\pe, x)\,\w{x_+^p}(\pe, x) - \w{\ln x_-}(\pe, x)\,\w{\de}(\pe, x)\right]\,dx = \frac{\sigma_p}{2} \,\lan\ps, \de\ran.\]
According to Definition 2, this proves equation (\ref{th1-}), and (\ref{th1+}) is obtained on the replacement $x \mapsto -x$.  The proof is complete.
 
\vspace*{1mm}
The above balanced products of the components $x_{\pm}^{-p}$ supported on the corresponding real half-lines can be employed further for obtaining results on singular products of the distribution $x^{-p}$. 

\begin{Th2} For each $p\in\NN$, the following balanced products hold in \GR \,:
\be \ \w{x^{\,-p-1}}\,.\,\w{x_{+}^p} \ + \ \w{\ln |x|}\,.\,\w{\de}(x) \ \approx \ x_{+}^{-1} \ - \ \si_p \ \de(x). \label{th2+}\ee
\be ( - 1)^{p+1} \ \w{x^{\,-p-1}}\,.\,\w{x_{-}^p} \ + \ \w{\ln |x|}\,.\,\w{\de}(x) \ \approx \ x_{-}^{-1} \ - \ \si_p \ \de(x). \label{th2-}\ee
\end{Th2}

\no \textsc{Proof}\,:   (i) Consider the following chain of identities and associations in \GR, taking into account equations (\ref{x+-p}) and (\ref{th1-})\,:
\bea \w{x_+^{\,-p-1}}\,.\,\w{x_+^p} & = & \w{x_+^{\,-p-1}}\,.\,\left( \w{x^{\,p}} + (-1)^p \
\w{x_-^{\,p}} \right) = - (-1)^p \ \w{x_+^{\,-p-1}}\,.\,\w{x_-^p} + \w{x_+^{\,-p-1}}\,.\,\w{x^{\,p}} \\ & \approx & - \ \w{\ln x_+}\,.\,\w{\de}(x) \ - \ \frac{\sigma_p}{2} \ \w{\de}(x) \ + \ \w{x_+^{-1}}. \eea
Here we have used $p$ times equation (\ref{xx}), as well as the fact that \ $f\,\w{u} \approx  f\,u$, \,for arbitrary $f\in C^{\infty}(\R)$ \,and $u\in \DD$ \cite[\S\,8.2]{chi}. \,The equivalence relation (\ref{equiv}) now yields
\be \w{x_+^{\,-p-1}}\,.\,\w{x_+^p} \ + \ \w{\ln x_+}\,.\,\w{\de}(x) \ \approx  \ x_+^{-1} \ - \ \frac{\si_p}{2} \ \de(x).  \label{55}\ee
(ii) Employing again equations (\ref{x+-p}) and (\ref{th1-}), as well as the balanced product (\ref{55}), we have
\bea \w{x^{\,-p-1}}\,.\,\w{x_+^p} & = & \left( \w{x_+^{\,-p-1}} + (- 1)^{p+1} \ \w{x_-^{\,-p-1}} \right).\,\w{x_+^p} \, = \,\w{x_+^{\,-p-1}}\,.\,\w{x_+^p} - (- 1)^p \w{x_-^{\,-p-1}}\,.\,\w{x_+^p}\\ &\approx & - \,\w{\ln x_+}\,.\,\w{\de}(x) - \w{\ln x_-}\,.\,\w{\de}(x)  +  \w{x_+^{-1}} - \si_p \,\w{\de}(x) = - \,\w{\ln |x|}\,.\,\w{\de}(x) + x_+^{-1} - \si_p \,\w{\de}(x).\eea
Whence,
\[ \w{x^{\,-p-1}}\,.\,\w{x_+^p} \ + \ \w{\ln |x|}\,.\,\w{\de}(x) \ \approx \w{x_+^{-1}} - \si_p \,\w{\de}(x).\]
In view of the equivalence relation (\ref{equiv}), this proves equation (\ref{th2+}); the replacement  $x \mapsto - x$ \ proves equation (\ref{th2-}) and the theorem. 

A direct consequence of Theorem 2 is the following. 

\begin{Cor1}$\!\!.$ For each $p\in\N$, the embeddings of the distributions $x^{-p-1}$ and $x_{\pm}^{\,p}$ satisfy\,:
\be \w{x^{\,-p-1}}\,.\,\w{x_{+}^p} \ - \ \w{x^{-1}}\,.\,\w{H} \ \approx \ - \ \si_p \  \de(x). \label{cor1+}\ee
\be ( - 1)^{p+1} \ \w{x^{\,-p-1}}\,.\,\w{x_{-}^p} \ - \ \w{x^{-1}}\,.\,\w{\check{H}} \ \approx \ \si_p\ \de(x). \label{cor1-}\ee
\end{Cor1}

Another implication from the result of Theorem 2 is given by this.

\begin{Cor2}$\!\!.$ The following balanced products hold for the embeddings in \GR \,of the distributions $(x {\pm}i 0)^{-p-1}$ \,and $x_+^{\,p}, \ p\in\N$\,: 
\be \w{(x {\pm}i 0)^{-p-1}}\,.\,\w{x_+^p} \ + \ \w{\ln |x|}\,.\,\w{\de}(x) \ \approx   \,x_+^{-1} \ - \left( \si_p \pm \,\frac{i \,\pi}{2} \right) \ \de(x). \label{cor2}\ee
\end{Cor2}

\no \textsc{Proof}\,: Employing equations (\ref{xip}), (\ref{th2+}), and the Colombeau product (\ref{cmuc}), one obtains
\bea \w{(x {\pm}i 0)^{-p-1}}\, .\,\w{x_+^p} \,& =& \, \w{x^{-p-1}}\, .\,\w{x_+^p} \ \mp \,\frac{( - 1)^p \ i \,\pi}{p!} \ \w{\de}^{(p)}(x)\, .\,\w{x_+^p} \\ &\ \approx & \,- \ \w{\ln |x|}\,.\,\w{\de}(x)\,+ \,\w{x_+^{-1}} \ - \ \si_p  \ \w{\de}(x) \ \mp \,i \,\pi\,/\,2 \ \w{\de}(x);\eea
which in view of the equivalence relation (\ref{equiv}) \,proves (\ref{cor2}). 

Consider further equation (\ref{th2+}) \,in the particular case $p=0$\,:
\be  \w{x^{-1}}\,.\,\w{H} \ + \ \w{\ln |x|}\,.\,\w{\de}(x) \ \approx \ x_{+}^{-1}. \label{th2-0}\ee
This equation will serve as a starting point for another generalization obtained by the next theorem. Its proof provides examples of balanced distributional products that are "stable under differentiation": the differentiation rule (\ref{l1}) leads again to balanced products (which is not true in general).

\begin{Th3}$\!\!.$ For each $p\in\NN$, the embeddings in \GR \,of the distributions $x^{-p-1}, H$, and $\de^{\,(p)}(x)$ \,satisfy 
\be \w{x^{\,-p-1}}.\,\w{H} + \frac{(-1)^p}{p!}\ \w{\ln |x|}\,.\,\w{\de}^{\,(p)}(x) \ \approx \ x_+^{-p-1}. \label{th3+}\ee  
\be (-1)^{p-1} \ \w{x^{\,-p-1}}\, .\,\w{\check{H}} \,+ \ \frac{1}{p!}\ \w{\ln |x|}\,.\,\w{\de}^{\,(p) }(x) \ \approx \ x_-^{-p-1}. \label{th3-}\ee  
\end{Th3}

\no \textsc{Proof}\,:  We shall make use of equation (\ref{xpdq}) written in the particular case $q=1$\,:
\be \w{x^{-p}}\,.\,\w{\de}(x) \ + \ \frac{(-1)^{p-1}}{(p-1)!} \ \w{x^{- 1}}\,.\,\w{\de}^{\,(p-1)}(x) \approx \frac{(-1)^{p}}{p!} \,\de^{(p)}(x), \qquad p\in\N. \label{xpd1}\ee

\no Apply now the differentiation rule (\ref{l1}) to the balanced product (\ref{th2-0}). Taking into account equations (\ref{defx-p}), as well as relation (\ref{cons}) for the consistency of differentiation with the embedding of distributions, we obtain
\[ - \,\w{x^{\,-2}}\, .\,\w{H} \,+ \,2\ \w{x^{- 1}}\,.\,\w{\de}(x) \,+ \,\w{\ln |x|}\,.\,\w{\de'}(x)\ \approx \ - \,x_+^{-2} \,- \,\de'(x).\]
On the strength of equation (\ref{xpd1}) in the case $p=1$, it follows
\be \w{x^{\,-2}}\, .\,\w{H} \,- \,\w{\ln |x|}\,.\,\w{\de'}(x)\ \approx \ x_+^{-2}.\label{p=1}\ee
Further differentiation of the latter equation according to (\ref{l1})  \,yields
\[ - \,2\,\w{x^{\,-3}}\, .\,\w{H} \,+ \,\w{x^{- 2}}\,.\,\w{\de}(x)\, -  \,\w{x^{- 1}}\,.\,\w{\de'}(x) + \,\w{\ln |x|}\,.\,\w{\de''}(x)\ \approx \ - 2\,x_+^{-3} \,+ \,\frac{1}{2}\,\de''(x).\]
Then equation (\ref{xpd1}) in the case $p=2$ \,allows us to replace the balanced product \,$\w{x^{- 2}}\,.\,\w{\de}(x)\, -  \,\w{x^{- 1}}\,.\,\w{\de'}(x)$ with the associated distribution $1/2\,\de''(x)$, which gives
\be \w{x^{\,-3}}\, .\,\w{H} \,+ \,\w{\ln |x|}\,.\,\w{\de''}(x)\ \approx \ x_+^{-3}.\label{p=2}\ee

\no This procedure can be repeated further, so we suppose the following balanced product holds that coincides with equations (\ref{p=1}), (\ref{p=2}) when $p=1, 2$\,:
\[ \w{x^{\,-p}}\, .\,\w{H} \,+ \,\frac{(-1)^{(p-1)}}{(p-1)!}\ \w{\ln |x|}\,.\,\w{\de^{\,(p-1)}}(x) \ \approx \ x_+^{-p}. \] 

\no Differentiation of this product according to (\ref{l1}) yields
\[ \w{x^{\,-p-1}} .\w{H} + \w{x^{-p}}.\w{\de}(x)  +  \frac{(-1)^{p-1}}{(p-1)!} \,\w{x^{- 1}}.\w{\de}^{\,(p)}(x) + \frac{(-1)^p}{p!}\ \w{\ln |x|}\,.\,\w{\de}^{\,(p)}(x) \,\approx \,x_+^{-p-1} +  \frac{(-1)^{p}}{p!}\,\de^{(p)}(x). \]
Applying then equation (\ref{xpd1}), we prove by induction equation (\ref{th3+}) for arbitrary $p\in\NN$. The replacement $x \mapsto -x$ in (\ref{th3+}) proves equation (\ref{th3-}) and the theorem.  

\vspace*{2mm}
\setlength{\baselineskip}{18pt}

\vspace*{2mm} 
\setlength{\baselineskip}{16pt}
\no {\footnotesize{Bulgarian Acad. Sci., INRNE - Theory Group, 72 Tzarigradsko shosse, 1784 Sofia, Bulgaria\\ E-mail: damyanov@netel.bg}}

\end{document}